\documentclass[a4paper,12pt]{article}
\usepackage{amsmath,amssymb,amsthm}
\usepackage{mathrsfs}
\newtheorem{theorem}{Theorem}[section]
\newtheorem{proposition}{Proposition}[section]
\newtheorem{lemma}{Lemma}[section]
\newtheorem{corollary}{Corollary}[section]
\newtheorem{definition}{Definition}[section]
\numberwithin{equation}{section}
\date{}
\title{\textbf{Some intrinsic properties of h-Randers conformal change}}
\begin{document}
\maketitle
\begin{center}
\textbf{\author{$ ^{1} $H. S. Shukla, $ ^{2} $V. K. Chaubey and $ ^{3} $Arunima Mishra}}
\end{center}
\begin{center}
$ ^{1} $Department of Mathematics and Statistics, D. D. U. Gorakhpur \\* University, Gorakhpur (U.P.)-273009,$ \; $India, \\* E-mail:- profhsshuklagkp@rediffmail.com  \\*
\end{center}
\begin{center}
$ ^{2} $Department of Applied Sciences, Buddha Institute of Technology,\\* Sector-7, Gida, Gorakhpur,(U.P.) India'\\* E-Mail: vkchaubey@outlook.com\\*
\end{center}
\begin{center}
$ ^{3} $Department of Mathematics, St. Joseph's College for Women,\\* Civil Lines, Gorakhpur, (U.P.), India, \\* E-Mail: arunima16oct@hotmail.com
\end{center}

\begin{abstract}
   In the present paper we have considered h-Randers conformal change of a Finsler metric $ L $, which is  defined as
\begin{center}
$ L(x,y)\rightarrow \bar{L}(x, y)=e^{\sigma(x)}L(x, y)+\beta (x, y), $
\end{center}
where $ \sigma(x) $ is a function of x, $\beta(x, y) = b_{i}(x, y)y^{i}$ is a 1- form on $M^{n}$ and $b_{i}$ satisfies the condition of being an h-vector. We have obtained the expressions for geodesic spray coefficients under this change. Further we have studied some special Finsler spaces namely quasi-C-reducible, C-reducible, S3-like and S4-like Finsler spaces arising from this metric. We have also obtained the condition under which this change of metric leads a Berwald (or a Landsberg) space into a space of the same kind.\newline\\*
\textbf{Mathematics subject Classification:} 53B40, 53C60.\newline\\*
\textbf{Keywords:} h-vector; special Finsler spaces; geodesic; conformal change. 
\end{abstract}

\section{Introduction}
\hspace{10pt} Let $M^{n}$ be an n-dimensional differentiable manifold and $F^{n}$ be a Finsler space equipped with a fundamental function $L(x, y), (y^{i} = \dot{x}^{i})$ of $M^{n}$. If a differential 1-form $\beta(x, y) = b_{i}(x)y^{i}$ is given on $M^{n}$, M. Matsumoto \cite{M1} introduced another Finsler space whose fundamental function is given by
\begin{center}
$ \bar{L}(x,y)= L(x, y)+\beta(x, y) $
\end{center}
This change of Finsler metric has been called $ \beta $-change \cite{PS, S}.\newline \\*
\hspace{10pt} The conformal theory of Finsler spaces was initiated by M.S. Knebelman \cite{Kne} in 1929 and has been investigated in detail by many authors \cite{Has, Iz1, Iz2, Kit}. The conformal change is defined as
\begin{center}
$ \bar{L}(x,y)\rightarrow e^{\sigma(x)}L(x, y) $,
\end{center}
where $ \sigma(x) $ is a function of position only and known as conformal factor.\newline \\*
In 1980, Izumi \cite{Iz2} introduced the h-vector $ b_{i} $ which is v-covariantly constant with respect to Cartan's connection $ C\Gamma $ (i.e. $ b_{i}|_{j}=0 $) and satisfies the relation $ LC^{h}_{ij}b_{h}=\rho h_{ij} $, where $ C^{h}_{ij} $ are components of (h)hv-torsion tensor and $ h_{ij} $ are components of angular metric tensor. Thus the h-vector is not only a function of coordinates $ x^{i} $, but it is also a function of directional arguments satisfying $ L\dot{\partial}_{j}b_{i}=\rho h_{ij} $.\newline\\*
In the paper \cite{Abed} S. H. Abed generalized the above two changes and have introduced another Finsler metric named as Conformal $  \beta $- change and further Gupta and Pandey \cite{GP} renamed it Randers conformal change and obtained various important result in the filed of Finsler spaces.
Recently we \cite{SCM} have generalized the metric given by S. H. Abed  with the help of h-vector and have introduced another Finsler metric which is defined as    
\begin{equation}
\bar{L}(x, y)=e^{\sigma(x)}L(x, y)+\beta (x, y),
\end{equation}
where $ \sigma(x) $ is a function of x and $\beta(x, y) = b_{i}(x, y)y^{i}$ is a 1- form on $M^{n}$ and $ b_{i} $ satisfies the condition of being an h-vector, We call the change $ L(x, y)\rightarrow \bar{L}(x, y) $ as h-Randers conformal change. This change generalizes various types of changes. When $ \beta $ = 0, it reduces to a conformal change. When $ \sigma $ = 0, it reduces to a h-Randers change \cite{PPT}. When $ \beta=0 $ and $ \sigma $ is a non-zero constant then it reduces to a homothetic change. When $ b_{i} $ is function of position only and $ \sigma $ = 0, it reduces to Randers change\cite{SSAY}. When $ b_{i} $ and $ \sigma $ are functions of position only, it reduces to Randers conformal change \cite{Abed, GP}. \newline \\*
In the present paper we have obtained the expressions for geodesic spray coefficients under this change. Further we have studied some special Finsler spaces namely quasi C-reducible, C-reducible, S3-like and S4-like Finsler spaces arising from this metric. We have also obtained the conditions under which this change of metric leads a Berwald (or a Landsberg) space into a space of the same kind.
\section{h-Randers conformal change}
Let the Cartan's connection of Finsler space $ F^{n} $ be denoted by $ C\Gamma = (F^{i}_{jk}, G^{i}_{j}, C^{i}_{jk}) $. Since $ b_{i}(x, y) $ are components  of h-vector, we have
\begin{equation}
(a)\;\; b_{i}|_{j}=\dot{\partial}_{j}b_{i}-b_{h}C^{h}_{ij}=0\;\;\;\;\;\;\; (b)\;\; LC^{h}_{ij}b_{h}=\rho h_{ij}
\end{equation}
Hence we obtain
\begin{equation}
\dot{\partial}_{j}b_{i}=L^{-1}\rho h_{ij}
\end{equation}
Since $ h_{ij} $ are components of an indicatory tensor i.e. $ h_{ij}y^{j}=0 $, we have $ \dot{\partial}_{i}\beta=b_{i} $.
\begin{definition}
Let $M^{n}$ be an n-dimensional differentiable manifold and $F^{n}$ be a Finsler space equipped with a fundamental function $L(x, y), (y^{i} = \dot{x}^{i})$ of $M^{n}$. A change in the fundamental function $ L $ by the equation (1.1) on the same manifold $M^{n}$ is called h-Randers conformal change. A space equipped with fundamental metric $ \bar{L} $ is called  h-Randers conformally changed Finsler space $ \bar{F}^{n} $.  
\end{definition}
Differentiating equation (1.1) with respect to $ y^{i} $, the normalized supporting element $ \bar{l}_{i}=\dot{\partial}_{i}\bar{L} $ is given by
\begin{equation}
\bar{l}_{i}=e^{\sigma}l_{i}+b_{i},
\end{equation}
where $ l_{i}=\dot{\partial}_{i}L $ is the normalized supporting element $ l_{i} $ of $ F^{n} $.\newline \\*
Differentiating (2.3) with respect to $ y^{j} $ and using (2.2) and the fact that $ \dot{\partial}_{j}l_{i}=L^{-1}h_{ij} $, we get
\begin{equation}
\bar{h}_{ij}=\phi h_{ij},
\end{equation}
where $ \phi=L^{-1}\bar{L}(e^{\sigma}+\rho) $ and $ h_{ij}= L\dot{\partial}_{i}\dot{\partial}_{j}L $ is the angular metric tensor in the Finsler space $ F^{n} $. \newline \\*
Since $ h_{ij}=g_{ij}-l_{i}l_{j} $, from (2.3) and (2.4) the fundamental tensor $ \bar{g}_{ij}=\dot{\partial}_{i}\dot{\partial}_{j}\frac{\bar{L}^{2}}{2}=\bar{h}_{ij} +\bar{l}_{i}\bar{l}_{j} $ is given as
\begin{equation}
\bar{g}_{ij}=\phi g_{ij}+b_{i}b_{j}+e^{\sigma}(b_{i}l_{j}+b_{j}l_{i})+( e^{2\sigma}-\phi)l_{i}l_{j}
\end{equation}
It is easy to see that the det$( \bar{g}_{ij})$ does not vanish, and the reciprocal tensor with components $\bar{g}^{ij}$ of $ \bar{F}^{n} $, obtainable from $ \bar{g}^{ij}\bar{g}_{jk}=\delta^{i}_{k} $, is given by
\begin{equation}
\bar{g}^{ij}=\phi^{-1}g^{ij}-\mu l^{i}l^{j}-\phi^{-2}(e^{\sigma}+\rho)(l^{i}b^{j}+l^{j}b^{i}),
\end{equation}
where $ \mu =(e^{\sigma}+\rho)^{2}\phi^{-3}(e^{\sigma}-b^{2}-\phi) $, $ b^{2}=b_{i}b^{i} $, $ b^{i}=g^{ij}b_{j} $ and $ g^{ij} $ is the reciprocal tensor of $ g_{ij} $ of $ F^{n} $.\newline \\*
We have following lemma \cite{SCM}:
\begin{lemma}
The scalar $ \rho $ used in the condition of h-vector is a function of coordinates $ x^{i} $ only.
\end{lemma}
From equations (1.1), (2.3) and lemma 2.1 we have
\begin{equation}
\dot{\partial}_{i}\phi=L^{-1}(e^{\sigma}+\rho)m_{i},
\end{equation}
where
\begin{equation}
m_{i}=b_{i}-(L^{-1}\beta)l_{i}
\end{equation}
Differentiating (2.4) with respect to $ y^{k} $ and using (2.3), (2.4), (2.7) and the relation $ \dot{\partial}_{k}h_{ij}=2C_{ijk}-L^{-1}(l_{i}h_{jk}+l_{j}h_{ik}) $, the Cartan covariant tensor $\bar{C}_{ijk}$ is given by
\begin{equation}
\bar{C}_{ijk}=\phi C_{ijk}+\frac{(e^{\sigma}+\rho)}{2L}(h_{ij}m_{k}+h_{jk}m_{i}+ h_{ki}m_{j}),
\end{equation}
where $ C_{ijk} $ is (h)hv-torsion tensor of Cartan's connection $ C\Gamma $ of Finsler space $ F^{n} $.\newline \\*
From the definition of $ m_{i} $, it is evident that
\begin{eqnarray}
(a)\;\; m_{i}l^{i}=0, \;\;\;\;\; (b)\;\;m_{i}b^{i}=b^{2}-\frac{\beta^{2}}{L^{2}}=b_{i}m^{i},\;\;\;\;\;\;\; \\\nonumber
(c)\;\;g_{ij}m^{i}=h_{ij}m^{i}=m_{j},\;\;\;\;\;(d)\;\;C_{ihj}m^{h}=L^{-1}\rho h_{ij} 
\end{eqnarray}
From (2.1), (2.6), (2.9) and (2.10), we get 
\begin{eqnarray}
\bar{C}^{h}_{ij}=C^{h}_{ij}+\frac{1}{2\bar{L}}(h_{ij}m^{h}+h^{h}_{j}m_{i}+ h^{h}_{i}m_{j})-\frac{1}{\bar{L}}[\lbrace \rho +\frac{L}{2\bar{L}}(b^{2} -\\\nonumber \frac{\beta^{2}}{L^{2}})\rbrace h_{ij}+\frac{L}{\bar{L}}m_{i}m_{j}]l^{h}
\end{eqnarray} 
\begin{proposition}
Let $ \bar{F}^{n} = (M^{n}, \bar{L})$ be an n-dimensional Finsler space obtained from the h-Randers conformal change of the Finsler space $F^{n} = (M^{n}, L)$, then the normalized supporting element $\bar{l}_{i}$, angular metric tensor $ \bar{h}_{ij}$, fundamental metric tensor $\bar{g}_{ij}$  and (h)hv-torsion tensor $\bar{C}_{ijk}$ of $ \bar{F}^{n} $ are given by (2.3), (2.4), (2.5) and (2.9) respectively.
\end{proposition}

\section{Geodesic Spray coefficients of $ \bar{F}^{n} $}
\hspace{10pt} Let $ s $ be the arc-length of a curve $ x^{i}=x^{i}(t) $ on a differentiable manifold $ M^{n} $, then the equation of a geodesic \cite{KAM} of $ F^{n}=(M^{n}, L) $ is written in the well-known form:
\begin{equation}
\frac{d^{2}x^{i}}{ds^{2}}+2G^{i}(x, \frac{dx}{ds})=0,
\end{equation}
where functions $ G^{i}(x, y) $ are the geodesic spray coefficients given by
\begin{center}
$ 2G^{i}=g^{ir}(y^{j}\dot{\partial}_{r}\partial_{j}F-\partial_{r}F), \;\;\;\;\;\;\;\; F=\frac{L^{2}}{2}. $
\end{center}
Now suppose $ \bar{s} $ is the arc-length of a curve  $ \bar{x}^{i}=\bar{x}^{i}(t) $ on a differentiable manifold $ M^{n} $ in the Finsler space $ \bar{F}^{n}=(M^{n}, \bar{L}) $, then the equation of geodesic in $ \bar{F}^{n} $ can be written as
\begin{equation}
\frac{d^{2}x^{i}}{d\bar{s}^{2}}+2\bar{G}^{i}(x, \frac{dx}{d\bar{s}})=0,
\end{equation}
where functions $ \bar{G}^{i}(x, y) $ are given by
\begin{center}
$ 2\bar{G}^{i}=\bar{g}^{ir}(y^{j}\dot{\partial}_{r}\partial_{j}\bar{F}-\partial_{r}\bar{F}), \;\;\;\;\;\;\;\; \bar{F}=\frac{\bar{L}^{2}}{2}. $
\end{center}
Since $ d\bar{s}=\bar{L}(x, dx) $, this is also written as
\begin{center}
$ d\bar{s}=e^{\sigma(x)}L(x, dx)+b_{i}(x, y)dx^{i}=e^{\sigma(x)}ds+b_{i}(x, y)dx^{i} $
\end{center}
Since $ ds=L(x, dx) $, we have
\begin{equation}
\frac{dx^{i}}{ds}=\frac{dx^{i}}{d\bar{s}}[e^{\sigma(x)}+b_{i}(x, y)\frac{dx^{i}}{ds}]
\end{equation}
Differentiating (3.3) with respect to $ s $, we have
\begin{center}
$ \frac{d^{2}x^{i}}{ds^{2}}=\frac{d^{2}x^{i}}{d\bar{s}^{2}}[e^{\sigma(x)}+b_{i}\frac{dx^{i}}{ds}]^{2}+\frac{dx^{i}}{d\bar{s}}(\frac{de^{\sigma(x)}}{ds}+\frac{db_{i}}{ds}\frac{dx^{i}}{ds}+b_{i}\frac{d^{2}x^{i}}{ds^{2}}) $
\end{center}
Substituting the value of $ \frac{dx^{i}}{d\bar{s}} $ from (3.3), the above equation becomes
\begin{eqnarray}
\frac{d^{2}x^{i}}{ds^{2}}=\frac{d^{2}x^{i}}{d\bar{s}^{2}}[e^{\sigma(x)}+b_{i}\frac{dx^{i}}{ds}]^{2}+\frac{\frac{dx^{i}}{ds}}{[e^{\sigma(x)}+b_{i}\frac{dx^{i}}{ds}]}(\frac{de^{\sigma(x)}}{ds}+\\\nonumber \frac{db_{i}}{ds}\frac{dx^{i}}{ds}+b_{i}\frac{d^{2}x^{i}}{ds^{2}})
\end{eqnarray}
Now differentiating equation (1.1) with respect to $ x^{i} $ we have
\begin{equation}
\partial_{i}\bar{L}=e^{\sigma}A_{i}+B_{i},
\end{equation}
where $ A_{i}=L\partial_{i}\sigma + \partial_{i}L $ and $ B_{i}=\partial_{i}\lbrace b_{r}(x,y)\rbrace y^{r} $.\newline \\*
Differentiating above equation with respect to $ y^{j} $ we have
\begin{equation}
\dot{\partial}_{j}\partial_{i}\bar{L}=e^{\sigma}\dot{\partial}_{j}A_{i} + \dot{\partial}_{j}B_{i},
\end{equation}
where $ \dot{\partial}_{j}A_{i}= l_{j}\partial_{i}\sigma +\dot{\partial}_{j}\partial_{i}L $ and $ \dot{\partial}_{j}B_{i}= \dot{\partial}_{j}\lbrace \partial_{i}b_{r}(x, y)\rbrace y^{r} + \partial_{i} b_{r}(x, y)\delta^{r}_{j} $. \newline \\*
Since 
\begin{center}
$ 2\bar{G}_{r}=y^{j}(\bar{l}_{r}\partial_{j}\bar{L}+ \bar{L}\dot{\partial}_{r}\partial_{j}\bar{L})-\bar{L}\partial_{r}\bar{L} $
\end{center}
therefore using equations (2.3), (3.5) and (3.6) we have
\begin{eqnarray}
2\bar{G}_{r}=2e^{2\sigma}G_{r}+y^{j}\lbrace 2e^{2\sigma}l_{r}L\partial_{j}\sigma +e^{\sigma}(l_{r}B_{j}+b_{r}A_{j})+b_{r}B_{j}+\\\nonumber e^{\sigma}L\dot{\partial}_{r}B_{j} +\beta e^{\sigma}\dot{\partial}_{r}A_{j}+\beta \dot{\partial}_{r}B_{j}\rbrace - (e^{2\sigma}L^{2}\partial_{r}\sigma + e^{\sigma}LB_{r}\\\nonumber +\beta e^{\sigma}A_{r}+\beta B_{r}),
\end{eqnarray}
where $ G_{r}=y^{j}\lbrace l_{r}\partial_{j}L+L\dot{\partial}_{r}\partial_{j}L \rbrace -L\partial_{r}L $ is the spray coefficients for the Finsler space $ F^{n} $.\newline \\*
Using equations (2.6) and (3.7) we have
\begin{equation}
\bar{G}^{i}=JG^{i}+M^{i},
\end{equation}
where $G^{i}=g^{ir}G_{r},\;\;\; J=\frac{1}{\phi},\;\;\;and\;\; M^{i}=\frac{1}{2}e^{2\sigma}G_{r}\lbrace -\mu l^{i}l^{r}-\phi^{-2}(e^{\sigma}+\rho)(l^{i}b^{r}+l^{r}b^{i}) \rbrace +\frac{1}{2}[\phi^{-1}g^{ir}-\mu l^{i}l^{r}-\phi^{-2}(e^{\sigma}+\rho)(l^{i}b^{r}+l^{r}b^{i})][y^{j}\lbrace 2e^{2\sigma}l_{r}L\partial_{j}\sigma +e^{\sigma}(l_{r}B_{j}+b_{r}A_{j})+b_{r}B_{j}+ e^{\sigma}L\dot{\partial}_{r}B_{j} +\beta e^{\sigma}\dot{\partial}_{r}A_{j}+\beta \dot{\partial}_{r}B_{j}\rbrace - (e^{2\sigma}L^{2}\partial_{r}\sigma +e^{\sigma}LB_{r} +\beta e^{\sigma}A_{r}+\beta B_{r})]  $.
\begin{theorem}
Let $\bar{F}^{n} = (M^{n}, \bar{L})$ be an n-dimensional Finsler space obtained from the h-Randers conformal change of the Finsler space $F^{n} = (M^{n}, L)$, then the the geodesic spray coefficients $ \bar{G}^{i} $ for the Finsler space $ \bar{F}^{n} $ are given by (3.8) in the terms of the geodesic spray coefficients $ G^{i} $ of the Finsler space $F^{n}$.
\end{theorem}
\begin{corollary}
Let $\bar{F}^{n} = (M^{n}, \bar{L})$ be an n-dimensional Finsler space obtained from the h-Randers conformal change of the Finsler space $F^{n} = (M^{n}, L)$, then the equation of geodesic of $ \bar{F}^{n} $ is given by (3.2), where $ \frac{d^{2}x^{i}}{d\bar{s}^{2}} $ and $ \bar{G}^{i} $ are given by (3.4) and (3.8) respectively.
\end{corollary}
\section{C-reducibilty of $ \bar{F}^{n} $}
\hspace{10pt} Following Matsumoto \cite{M2}, in this section we shall investigate special cases of the Finsler space with h-Randers conformally changed Finsler space $ \bar{F}^{n} $.
\begin{definition}
A Finsler space $(M^{n}, L)$ with dimension $n \geq 3$ is said to be quasi-C-reducible if the Cartan tensor $C_{ijk}$ satisfies
\begin{equation}
C_{ijk} = Q_{ij}C_{k} + Q_{jk}C_{i} + Q_{ki}C_{j},
\end{equation}
where $Q_{ij}$ is a symmetric indicatory tensor.
\end{definition}
Substituting $ h = j $ in equation (2.11) we get
\begin{equation}
\bar{C}_{i}=C_{i}+\frac{(n+1)}{2\bar{L}}m_{i}
\end{equation}
Using equations (2.9) and (4.2), we have
\begin{center}
$ \bar{C}_{ijk}=\phi C_{ijk}+\frac{\phi}{(n+1)}\pi_{(ijk)}\lbrace h_{ij}(\bar{C}_{k}-C_{k})\rbrace $,
\end{center}
where $ \pi_{(ijk)} $ represents cyclic permutation and sum over the indices $ i, j $ and $ k $.\newline \\*
The above equation can be written as
\begin{center}
$ \bar{C}_{ijk}= \phi C_{ijk}+\frac{\phi}{(n+1)}\pi_{(ijk)}(h_{ij}\bar{C}_{k})-\frac{\phi}{(n+1)}\pi_{(ijk)}(h_{ij}C_{k})$
\end{center}
Thus
\begin{lemma}
In an h-Randers conformally changed Finsler space $ \bar{F}^{n} $, the Cartan's tensor can be written in the form
\begin{equation}
\bar{C}_{ijk}=\pi_{(ijk)}(\bar{H}_{ij}\bar{C}_{k})+V_{ijk},
\end{equation}
where $ \bar{H}_{ij}=\frac{\bar{h}_{ij}}{(n+1)} $ and $ V_{ijk}= \phi C_{ijk}-\frac{\phi}{(n+1)}\pi_{(ijk)}(h_{ij}C_{k}) $.
\end{lemma}
Since $ \bar{H}_{ij} $ is a symmetric and indicatory tensor, so from the above lemma and (4.1) we get
\begin{theorem}
An h-Randers conformally changed Finsler space $ \bar{F}^{n} $ is quasi-C-reducible if the tensor $ V_{ijk} $ of equation (4.3) vanishes identically.
\end{theorem} 
We obtain a generalized form of Matsumoto's result known \cite{M2} as a corollary of the above theorem
\begin{corollary}
If $ F^{n} $ is Reimannian then an h-Randers conformally changed Finsler space $ \bar{F}^{n} $ is always a quasi-C-reducible Finsler space.
\end{corollary} 
\begin{definition}
A Finsler space $(M^{n}, L)$ of dimension $n \geq 3$ is called C-reducible if the Cartan tensor $C_{ijk}$ is written in the form
\begin{equation}
C_{ijk} =\frac{1}{(n+1)}(h_{ij}C_{k} + h_{ki}C_{j} + h_{jk}C_{i})
\end{equation}
\end{definition}
Now from equation (2.9) and definition of C-reducibility we have
\begin{equation}
\phi C_{ijk}=\pi_{(ijk)}(\bar{h}_{ij}N_{k}),
\end{equation}
where $ N_{k}=\frac{1}{(n+1)}\bar{C}_{k}-\frac{1}{2\bar{L}}m_{k} $. Conversely, if (4.5) is satisfied for certain covariant vector $ N_{k} $ then from (2.9) we have
\begin{equation}
\bar{C}_{ijk}=\frac{1}{(n+1)}\pi_{(ijk)}(\bar{h}_{ij}\bar{C}_{k})
\end{equation} 
Thus we have
\begin{theorem}
An h-Randers conformally changed Finsler space $ \bar{F}^{n} $ is C-reducible iff equation (4.5) holds good.
\end{theorem}
\begin{corollary}
If the Finsler space  $ F^{n} $ is C-reducible Finsler space, then an h-Randers conformally changed Finsler space $ \bar{F}^{n} $ is always a C-reducible Finsler space.
\end{corollary} 
\section{Some Important tensors of $ \bar{F}^{n} $}
\hspace{10pt} The $ v $-curvature tensor \cite{M2} of Finsler space with fundamental function $ L $ is given by
\begin{center}
$ S_{hijk}= C_{ijr}C^{r}_{hk}-C_{ikr}C^{r}_{hj} $
\end{center}
Therefore the $ v $-curvature tensor of an h-Randers conformally changed Finsler space $ \bar{F}^{n} $ will be given by
\begin{equation}
\bar{S}_{hijk}= \bar{C}_{ijr}\bar{C}^{r}_{hk}-\bar{C}_{ikr}\bar{C}^{r}_{hj}
\end{equation}
From equations (2.9) and (2.11) we have
\begin{eqnarray}
\bar{C}_{ijr}\bar{C}^{r}_{hk}=\phi [C_{ijr}C^{r}_{hk}+(\frac{\rho}{L\bar{L}}-\frac{m^{2}}{4\bar{L}^{2}})h_{hk}h_{ij}+\frac{1}{2\bar{L}}(C_{ijk}m_{h}+ \\\nonumber C_{ijh}m_{k}+C_{ihk}m_{j}+C_{hjk}m_{i})+\frac{1}{4\bar{L}^{2}}(h_{hj}m_{i}m_{k}+\\\nonumber  h_{hi}m_{j}m_{k}+h_{jk}m_{i}m_{h} +h_{ik}m_{h}m_{j})],
\end{eqnarray}
where $ h_{jr}C^{r}_{hk}=C_{jhk}=h^{r}_{j}C_{rhk},\;\;\;\;\;\;\; m_{i}m^{i}=m^{2} $.\newline \\*
Using equations (5.1) and (5.2) we have
\begin{eqnarray}
\bar{S}_{hijk}=\phi [S_{hijk}+(\frac{\rho}{L\bar{L}}-\frac{m^{2}}{4\bar{L}^{2}})\lbrace h_{hk}h_{ij}- h_{hj}h_{ik}\rbrace + \frac{1}{4\bar{L}^{2}}\lbrace h_{hj}m_{i}m_{k}\\\nonumber -h_{hk}m_{i}m_{j} + h_{ik}m_{h}m_{j}-h_{ij}m_{h}m_{k}\rbrace]
\end{eqnarray}
\begin{proposition}
In an h-Randers conformally changed Finsler space $ \bar{F}^{n} $ the v-curvature tensor $\bar{S}_{hijk}$ is given by (5.3).
\end{proposition}
It is well known\cite{M2} that the $ v $-curvature tensor of any three-dimensional Finsler space is of the form
\begin{equation}
L^{2}S_{hijk}=S(h_{hj}h_{ik}-h_{hk}h_{ij}),
\end{equation}
where scalar $ S $ in (5.4) is a function of $ x $ alone. \newline \\*
Owing to this fact M. Matsumoto defined the S3-like Finsler space as 
\begin{definition}
A Finsler space $ F^{n} $ $(n\geq 3)$ is said to be S3-like Finsler space if the $ v $-curvature tensor is of the form (5.4).
\end{definition} 
The $ v $-curvature tensor of any four-dimensional Finsler space may be written as \cite{M2}:
\begin{equation}
L^{2}S_{hijk}=\Theta_{(jk)}\lbrace h_{hj}K_{ki}+h_{ik}K_{hj} \rbrace ,
\end{equation}
where $ K_{ij} $ is a (0, 2) type symmetric Finsler tensor field which is such that $ K_{ij}y^{j}=0 $ and the symbol $ \Theta_{(jk)}\lbrace ... \rbrace $ denotes the interchange of $ j, k $ and subtraction. The definition of S4-like Finsler space is given as
\begin{definition}
A Finsler space $ F^{n}(n\geq 4) $ is said to be S4-like Finsler space if the $ v $-curvature tensor is of the form (5.5).
\end{definition} 
From equation (5.3) we have
\begin{lemma}
The v-curvature tensor $ \bar{S}_{hijk} $ of a h-Randers conformally changed Finsler space can be written as
\begin{equation}
\bar{S}_{hijk}=\bar{S}(\bar{h}_{hj}\bar{h}_{ik}-\bar{h}_{hk}\bar{h}_{ij})+ U_{hijk},
\end{equation}
where $ \bar{S}=-\frac{1}{\phi}(\frac{\rho}{L\bar{L}}-\frac{m^{2}}{4\bar{L}^{2}}) $ and $ U_{hijk}=\phi[S_{hijk}+\frac{1}{4\bar{L}^{2}}\lbrace h_{hj}m_{i}m_{k}-h_{hk}m_{i}m_{j} + h_{ik}m_{h}m_{j}-h_{ij}m_{h}m_{k}\rbrace] $
\end{lemma}
From lemma (5.1) and definition of S3-like Finsler space we have
\begin{theorem}
An h-Randers conformally changed Finsler space $ \bar{F}^{n} $ is S3-like if the tensor $ U_{hijk} $ of equation (5.6) vanishes identically. 
\end{theorem} 
From equation (5.3) we have 
\begin{lemma}
The v-curvature tensor $ \bar{S}_{hijk} $ of an h-Randers conformally changed Finsler space can also be written as
\begin{equation}
\bar{S}_{hijk}=\Theta_{(jk)}(\bar{h}_{ij}K_{ij}+\bar{h}_{ik}K_{hj}) + \phi S_{hijk},
\end{equation}
where $ K_{ij}=\frac{1}{4\bar{L}^{2}}m_{i}m_{j}-\frac{1}{2}(\frac{\rho}{L\bar{L}}-\frac{m^{2}}{4\bar{L}^{2}})h_{ij} $.
\end{lemma}
Thus from lemma(5.2) and definition of S4-like Finsler space we have
\begin{theorem}
If the $ v $-curvature tensor of Finsler space $ F^{n} $ vanishes identically then an h-Randers conformally changed Finsler space $ \bar{F}^{n} $ is S4-like Finsler space.
\end{theorem}
Now we are concerned with $ (v)hv $-torsion tensor $ P_{ijk} $. With respect to the Cartan's connection $ C\Gamma, \; L_{|i}=0,\;l_{i|j}=0,\; h_{ij|k}=0 $ hold good \cite{M2}.\newline \\*
Taking h-covariant derivative of the equation (2.9) we have
\begin{eqnarray}
\bar{C}_{ijk|h}=L^{-1}\bar{L}(e^{\sigma}\sigma_{|h}+\rho_{|h})\lbrace C_{ijk}+\frac{1}{2\bar{L}}(h_{ij}m_{k}+h_{jk}m_{i}+ h_{ki}m_{j})\rbrace \\\nonumber+ \phi \lbrace C_{ijk|h}+\frac{1}{2\bar{L}}(h_{ij}m_{k|h}+h_{jk}m_{i|h}+ h_{ki}m_{j|h})\rbrace ,
\end{eqnarray}
where $ m_{i|h}=b_{i|h}-L^{-1}l_{i}b_{r|h}y^{r} $.
\begin{lemma}
The $ h $-covariant derivative of the Cartan tensor $ \bar{C}_{ijk} $ of an h-Randers conformally changed Finsler space $ \bar{F}^{n} $ can be written as
\begin{equation}
\bar{C}_{ijk|h}=\phi C_{ijk|h}+V_{ijkh},
\end{equation}
where $ V_{ijkh}=L^{-1}\bar{L}(e^{\sigma}\sigma_{|h}+\rho_{|h})\lbrace C_{ijk}+\frac{1}{2\bar{L}}(h_{ij}m_{k}+h_{jk}m_{i}+ h_{ki}m_{j})\rbrace+\frac{\phi}{2\bar{L}}(h_{ij}m_{k|h}+h_{jk}m_{i|h}+ h_{ki}m_{j|h}). $
\end{lemma}
The $(v)hv$-torsion tensor $P_{ijk}$ of the Cartan connection $C\Gamma$ is written in the form 
\begin{center}
$ P_{ijk}=C_{ijk|0} $,
\end{center}
where the subscript '0' means the contraction with respect to the supporting element $ y^{i} $.\newline \\*
From the equation (5.8), the $(v)hv$-torsion tensor $\bar{P}_{ijk}$ is given by
\begin{eqnarray}
\bar{P}_{ijk}=\phi P_{ijk}+L^{-1}\bar{L}(e^{\sigma}\sigma_{|0}+\rho_{|0})\lbrace C_{ijk}+\frac{1}{2\bar{L}}(h_{ij}m_{k}+h_{jk}m_{i}+\\\nonumber h_{ki}m_{j})\rbrace + \frac{\phi}{2\bar{L}}\lbrace h_{ij}m_{k|0}+h_{jk}m_{i|0}+ h_{ki}m_{j|0}\rbrace 
\end{eqnarray}
Thus we have 
\begin{proposition}
The $(v)hv$-torsion tensor $\bar{P}_{ijk}$ of an h-Randers conformally changed Finsler space can be written in the form of (5.10).
\end{proposition}
From the equation (5.10) we have
\begin{lemma}
The $(v)hv$-torsion tensor $\bar{P}_{ijk}$ of an h-Randers conformally changed Finsler space can also be written as
\begin{equation}
\bar{P}_{ijk}=\phi P_{ijk}+W_{ijk},
\end{equation}
where $\; W_{ijk}=L^{-1}\bar{L}(e^{\sigma}\sigma_{|0}+\rho_{|0})\lbrace C_{ijk}+\frac{1}{2\bar{L}}(h_{ij}m_{k}+h_{jk}m_{i}+ h_{ki}m_{j})\rbrace + \frac{\phi}{2\bar{L}}\lbrace h_{ij}m_{k|0}+h_{jk}m_{i|0}+ h_{ki}m_{j|0}\rbrace. $  
\end{lemma}
We have
\begin{definition}
A Finsler space is called a Berwald space if $ C_{ijk|h}=0 $ holds good.
\end{definition}
\begin{definition}
A Finsler space is called a Landsberg space if $ P_{ijk}=0 $ holds good.
\end{definition}
In view of above definition (5.3) and the lemma (5.3) we have
\begin{theorem}
If a Finsler space $ F^{n} $ is a Berwald space and the tensor $ V_{ijkh} $ of equation (5.9) vanishes identically then an h-Randers conformally changed Finsler space $ \bar{F}^{n} $ is a Berwald space.
\end{theorem}
In view of above definition (5.4) and the lemma (5.4) we have
\begin{theorem}
If a Finsler space $ F^{n} $ is a Landsberg space and the tensor $ W_{ijk} $ of equation (5.11) vanishes identically then an h-Randers conformally changed Finsler space $ \bar{F}^{n} $ is a Landsberg space.
\end{theorem}
%
%
%
%
%


\begin{thebibliography}{15}

\bibitem{Has} 
	M. Hashiguchi, 
	{\em On conformal transformation of Finsler metric},
	J. Math. Kyoto Univ., 16, (1976), 25-50.
	
\bibitem{Iz1}
	H. Izumi,
	{\em Conformal transformations of Finsler spaces I},
	Tensor, N.S., 31, (1977), 33-41.
	
\bibitem{Iz2}
	H. Izumi, 
	{\em Conformal transformations of Finsler spaces II. An h-conformally flat Finsler space},
	Tensor, N.S., 33, (1980), 337-359.
	
\bibitem{Kit}
	M. Kitayama,
	{\em Geometry of transformations of Finsler metrics},
	Hokkaido University of Education, Kushiro Campus, Japan 2000.
	
\bibitem{KAM} 
	M. Kitayama, M. Azuma, M. Matsumoto, 
	{\em On Finsler spaces with $ (\alpha, \beta) $-metric. Regularity, Geodesics and Main scalars},
	J. Hokkaido Univ., Edu., 46, 1, (1995), 1-10.
	
\bibitem{Kne} 
	M. S. Knebelman, 
	{\em Conformal geometry of generalized metric spaces},
	Proc. Nat. Acad. Sci.USA., 15, (1929), 33-41 and 376-379.

\bibitem{M1}
   M. Matsumoto,
    {\em On some transformation of locally Minkowskian space},
    Tensor, N. S. 22, (1971), 103-111.
	
\bibitem{M2}
	M. Matsumoto, 
	{\em Foundation of Finsler geometry and special Finsler spaces},
	Kaiseisha Press, Otsu, Japan 1986.
	
\bibitem{PPT} 
	B. N. Prasad, T. N. Pandey, D. Thakur,
	{\em H-Randers change of Finsler Metric},
	Investigations in Mathematical Sciences, 1, (2011), 71-84.

\bibitem{PS}
    B. N. Prasad, J. N. Singh,
    {\em Cubic transformation of Finsler spaces and n-fundamental forms of their hypersurfaces},
    Indian J.Pure Appl. Math, 20, 3 (1989), 242-249.

\bibitem{S}
	C. Shibata,
	{\em On invariant tensors of $ \beta $-changes of Finsler spaces}, 
	J. Math. Kyoto Univ., 24, (1984), 163-188.

\bibitem{SSAY}
	C. Shibata, H. Shimada, M. Azuma, H. Yasuda, 
	{\em On Finsler spaces with Randers metric}, 
	Tensor, N. S., 31, (1977), 219-226.
	
\bibitem{SCM} 
	H. S. Shukla, V. K. Chaubey, Arunima Mishra,
	{\em On Finsler spaces with h-Randers conformal change},
	Tensor, N. S., 74, (2013), 135-144. 
	
\bibitem{Abed} 
	S. H. Abed, 
	{\em Conformal $  \beta $- changes in Finsler spaces}, 
	Proc. Math. Phys. Soc. Egypt, 86 (2008),79-89, ArXivNo.: math. DG/060240.
	
\bibitem{GP} 
	M. K. Gupta, P. N. Pandey, 
	{\em On hypersurface of a Finsler space with Randers conformal metric}, 
	Tensor, N. S.,70(2008): 229-240.	

\end{thebibliography}
\end{document}